\documentclass[12pt,reqno]{amsart}

\usepackage{amssymb}

\theoremstyle{plain}
\newtheorem{theorem}{Theorem}[section]
\newtheorem{lemma}[theorem]{Lemma}

\theoremstyle{remark}
\newtheorem{rem}[theorem]{Remark}

\newcommand{\thmref}[1]{Theorem~\ref{#1}}

\newcommand\iv{^{-1}} 
\newcommand\sbl[1]{\langle#1\rangle}   
\newcommand\cL{\mathcal{L}}
\newcommand\Mgl{\mathcal{M}_{\ell}}
\newcommand\Del{\mathcal{D}}
\newcommand\qdP{\times_Q}

\newcommand\SL{\mathrm{SL}}
\newcommand\PSL{\mathrm{PSL}}
\newcommand\SO{\mathrm{SO}}
\newcommand\SU{\mathrm{SU}}
\newcommand\PSp{\mathrm{PSp}}

\newcommand\RRR{\mathbb{R}}
\newcommand\CCC{\mathbb{C}}
\newcommand\Hq{\mathbb{H}}
\newcommand\homega{\widehat{\omega}}
\newcommand\conj[2]{#1#2#1^{-1}}

\title[Infinite Simple Bol Loops]
{Infinite Simple Bol Loops}

\author[H.~Kiechle]{Hubert Kiechle}
\address{Fachbereich Mathematik \\
Universit\"{a}t Hamburg \\
Bundesstrasse 55 \\
D-20146 Hamburg, Germany}
\email{kiechle@math.uni-hamburg.de}
\urladdr{http://www.math.uni-hamburg.de/home/kiechle/kiechle.html}
\author[M.~K.~Kinyon]{Michael~K.~Kinyon}
\address{Department of Mathematical Sciences \\
Indiana University South Bend \\
South Bend, IN 46634 USA}
\email{mkinyon@iusb.edu}
\urladdr{http://mypage.iusb.edu/\symbol{126}mkinyon}

\date{\today}

\begin{document}

\begin{abstract}
If the left multiplication group of a loop is simple, then the loop
is simple. We use this observation to give examples of infinite
simple Bol loops.
\end{abstract}

\maketitle

\section* {Introduction}

A loop satisfying the identity $x(y\cdot xz) = (x\cdot yx)z$ is a
(left) \emph{Bol loop}. It is generally agreed that the most significant
open problem in loop theory today is the existence or nonexistence of
finite, simple Bol loops which are not Moufang. In this note, we give examples
of infinite simple Bol loops. The construction of the loops is well-known, and
all we actually need to do is to compute the left multiplication group, which
turns out to be simple. It follows from a lemma in {\S}1 that the loop is simple,
as well.

In {\S}2 we recall the construction based on $\SL(n,F)$ over
suitable fields $F$ to obtain examples. In {\S}3 we use a result
of R\'ozga \cite{Rozga00} on Lie groups of noncompact type to obtain
even more examples.

The example from $\SL(2,\RRR)$ has been treated in detail in 
\cite[{\S}22]{NagyStr-Loops}. In particular, it is shown that
the topological left multiplication group is $\PSL(2,\RRR)$, a simple
group. Our results imply that $\PSL(2,\RRR)$ is also the (combinatorial)
left multiplication group, which is what we need for our conclusion.

It turns out that all of the examples we provide satisfy the
automorphic inverse property $(xy)\iv = x\iv y\iv$. Bol loops
satisfying this property are known by various names:
K-loops \cite{K-KL}, Bruck loops \cite{AKP}, or
gyrocommutative gyrogroups \cite{FU}. We do not know of any
simple, proper Bol loops not satisfying the automorphic inverse property.

In the finite case, it is known that a simple K-loop, if such a loop exists, must
consist only of elements of order a power of $2$ \cite{AKP}. By contrast, the
infinite simple K-loops described here are uniquely $2$-divisible, that is, the
squaring map is a permutation.

Our notation conventions follow \cite{K-KL}. Mappings act on the left
of their arguments. For a loop $L$, we
denote left translations by $\lambda_x y := x y$. The left multiplication
group is $\Mgl(L) := \sbl{\lambda_x : x\in L}$, while the left inner mapping
group is $\Del(L) := \Mgl(L)_1 = \{ \alpha \in \Mgl : \alpha(1) = 1 \}$.
$\Del(L)$ is generated by all mappings
$\delta_{x,y} := \lambda_{xy}\iv \lambda_x \lambda_y$, $x,y\in L$.
The quasidirect product of $L$ and $\Del(L)$, denoted by $L\qdP \Del(L)$,
is the set $L\times \Del(L)$ with the product $(x,\alpha)(y,\beta) :=
(x\cdot \alpha(y), \delta_{x,\alpha(y)}\alpha \beta )$. $L\qdP \Del(L)$ 
is just an isomorphic copy of $\Mgl$, but is more convenient for some
purposes.

\section {Left Multiplication Group}
\label{sec:LMG}

Let $L$ be a loop. A subloop $N$ of $L$ is called \emph{normal} if 
\[
ab\cdot N=a\cdot bN=aN\cdot b \text{\ for all\ }a,b\in L
\]
Then there is a \emph{factor loop} $L/N$ and a canonical epimorphism
$L\to L/N$ with kernel $N.$ Conversely, the kernel of any homomorphism
is a normal subloop, see \cite[IV]{B-BS}.
For each subloop $N$ in $L$, define
\[
\cL(N) := \{ \alpha\in\Mgl : \alpha(x)\in N x, \forall x\in L \}.
\]

As in \cite[IV.1]{B-BS} one shows

\begin{lemma}
\label{lem:basic}
Let $L,L'$ be loops and $\theta:L\to L'$ a homomorphism. Then
\begin{enumerate}
\item The map  $\theta_*:\Mgl(L)\to\Mgl(\theta(L))$, defined by
$\theta_*(\alpha)(\theta(x))=\theta\alpha(x)$ for all $x\in L$
is an epimorphism with kernel  $\cL(\ker\theta).$

\item  If $N$ is a normal subloop in $L,$ then $\cL(N)$
is a normal subgroup in $\Mgl(L)$ and $\Mgl(L/N)=\Mgl(L)/\cL(N).$
Moreover, $N$ is a proper subloop of
$L$ if and only if $\cL(N)$ is a proper subgroup of $\Mgl(L)$.
\end{enumerate}
\end{lemma}

Indeed, one computes $\theta_*(\lambda_a)=\lambda_{\theta(a)}$, and
the rest is straightforward. See also \cite[1.7]{NagyStr-Loops}.

A loop is called \emph{simple} if every normal subloop is trivial.
As an immediate corollary we get

\begin{theorem}
\label{thm:simple}
If $\Mgl$ is simple, then $L$ is simple.
\end{theorem}

\section {Examples}
\label{sec:examples}

We now will now use \thmref{thm:simple} to show the simplicity
of certain K-loops.

Let $R$ be an ordered field, and let $K:=R(i)$, where $i^2=-1.$
Then $K$ has a unique nonidentity automorphism
$z\mapsto \bar{z}$, which fixes $R$ elementwise and maps $i$ to $-i$.
This will also be applied componentwise to matrices over $K$.
Let $n\geq 2$ be a fixed integer and consider the set $L$ of positive
definite symmetric, resp. hermitian $n\times n$-matrices over $R$, resp., $K$.
We assume that $R$ is \emph{$n$-real}, i.e., the characteristic
polynomial of every matrix in $L$ splits over $K$ (in fact over $R$)
into linear factors. 

\begin{rem}
  An ordered field is 2-real if and only if it is
  pythagorean. Therefore, every $n$-real field is pythagorean. Real
  closed fields are $n$-real for every~$n$. See \cite{K-KL} for more
  on $n$-real fields.
\end{rem}

Furthermore, let $G=\SL(n,F),$ and $\Omega=\SO(n,R)$, or
$\Omega=\SU(n,K)$, where $F = R$, $K$, respectively.
Finally, put $L_G := L\cap G$.

Then $\Omega$ is a subgroup of $G$ and $L_G$ is a transversal containing
$I_n$, the $n\times n$ identity matrix. In particular, $G=L_G \Omega$.
Then for all $A,B\in L_G$ there exist unique $A\circ B\in L_G$, $d_{A,B}\in\Omega$
such that $AB=(A\circ B)d_{A,B}$.

By \cite[(9.2)]{K-KL} we have

\begin{theorem}
\label{thm:K-loop}
$(L_G,\circ)$ is a K-loop, and therefore a Bol loop.
\end{theorem}

Our aim is to show

\begin{theorem}
\label{thm:ex-simple}
$(L_G,\circ)$ is simple.
\end{theorem}

\begin{proof}
Indeed, we will show that $\Mgl = \PSL(n,F)$, which is a simple
group. The desired result will then follow from \thmref{thm:simple}.

Consider the map
\[
\Phi: \left \{
\begin{array}{rl}
 G=L_G\Omega &\to\; L_G\qdP \Del \\
 g=A\omega &\mapsto\; (A,\homega)
\end{array} \right.\ (A\in L_G,\,\omega\in\Omega),
\]
where $\homega(A)=\conj\omega A$ for all $A\in L_G$. Recall that
$\homega$ acts on $L_G$, that is, $\homega(A)\in L_G$ for all $A\in L_G$.
By \cite[(9.3.1)]{K-KL} we have $\Del=\{\homega;\,\omega\in\Omega\},$
so $\Phi$ is well-defined and surjective.

Let $g=A\omega,g'=B\omega'\in G$ be decomposed as in the definition of
$\Phi.$ Then
\[
gg' = A\omega B\omega' = A\homega(B) \omega\omega'
= (A\circ \homega(B)) d_{A,\homega(B)}\omega\omega' .
\]
Since $\widehat{d_{A,\homega(B)}} = \delta_{A,\homega(B)}$, we
have 
\[
\Phi(gg')=
\left(A\circ\homega(B),
   \delta_{A,\homega(B)}\homega\homega'\right) 
= (A,\homega)(B,\homega')=\Phi(g)\Phi(g'). 
\]
Thus $\Phi$ is an epimorphism. Then $N := \mathrm{ker}(\Phi) = 
\{ \omega \in \Omega : \widehat{\omega}(A) = A, \ \forall A\in L_G \}$,
that is, $N$ is the centralizer of $L_G$ inside $\Omega$. By \cite[(1.22)]{K-KL}
$N$ is the center of $\SL(n,F).$ Therefore,
\[
\Mgl = L_G \qdP \Del=G/N=\PSL(n,F)
\]
is a simple group (see, e.g., \cite[10.8.4]{S-GT} or  \cite[II.6.13]{HuppI}).
\end{proof}

\section{Examples from Simple Lie Groups}
\label{ex:LieGroups}

Let $G$ be a noncompact Lie group which is simple as an abstract
group. Then $G$ has a simple Lie algebra. By \cite[{\S}6]{Rozga00}, 
there exists a subgroup $\Omega < G$ and a subset $L\subset G$ such
that $L$ is a transversal of $G/\Omega$ and $L$ is a K-loop. Indeed,
$L$ and $\Omega$ come from the Cartan decomposition of the Lie algebra
of $G$.

\begin{theorem}
$L$ is a simple K-loop.
\end{theorem}

\begin{proof}
By \cite[Prop.~24]{Rozga00}, $L$ generates $G$.
Then by \cite[(2.8)]{K-KL}, $\Mgl(L)$ is a homomorphic image of
$G$. Since $G$ is simple, we have $\Mgl(L) = G$. An application of
\thmref{thm:simple} gives the result.
\end{proof}

This applies to the groups $\PSL(n,\RRR)$, $\PSL(n,\CCC)$, $\PSL(n,\Hq)$
as well as $\PSp(2n,\RRR)$, $\PSp(2n,\CCC)$. Simplicity for these groups is
proved for instance in \cite[II.6.13, II.9.22] {HuppI}.


\begin{thebibliography}{99}

\bibitem{AKP} M.~Aschbacher, M.K.~Kinyon, and J.D.~Phillips,
Finite Bruck loops, submitted. Available at
{\tt http://www.arXiv.org/abs/math.GR/0401193}.

\bibitem{B-BS} R.H.~Bruck,
\textit{A Survey of Binary Systems},
Springer-Verlag, Berlin-Heidelberg-New York, 1971.
MR 20{\#}76, Zbl. 206:30301.

\bibitem{FU} T.~Foguel and A.A.~Ungar,
Gyrogroups and the decomposition of groups into twisted subgroups and subgroups,
\textit{Pacific J. Math.} \textbf{197} (2001) 1--11.
MR 2002e:20142,  Zbl. pre01589578.

\bibitem{HuppI} B.~Huppert,
\textit{Endliche Gruppen} I,
Springer-Verlag, Berlin-Heidelberg-New York, 1967.
MR 37 {\#}302, Zbl. 0217.07201.

\bibitem{K-KL} H.~Kiechle,
\textit{Theory of {K}-loops},
Lecture Notes in Math. \textbf{1778},
Springer-Verlag, Berlin-Heidelberg-New York, 2002.
MR 2003d:20109, Zbl. 0997.20059

\bibitem{NagyStr-Loops} P.~T.~Nagy and K.~Strambach,
\textit{Loops in Group Theory and Lie Theory},
de Gruyter Expositions in Mathematics \textbf{35},
Walter de Gruyter, Berlin-New York, 2003.
MR 2003d:20110, Zbl. pre01732502.

\bibitem{Rozga00} K.~R{\'o}zga,
On central extensions of gyrocommutative gyrogroups,
\textit{Pacific J. Math.} \textbf{193} (2000) 201--218.
MR 2001a:20115,  Zbl. 1010.20055.

\bibitem{S-GT} W.~R.~Scott,
\textit{Group Theory},
Dover, New York, 1987.
MR 88d:20001, Zbl. 0641.20001.

\end{thebibliography}
\end{document}